\documentclass[12pt,a4paper]{article}
\usepackage[T2A]{fontenc}
\usepackage[cp1251]{inputenc}
\usepackage[english]{babel}
\usepackage{epsfig,amssymb,amsmath}

\begin{document}

\title{The wild Fox-Artin arc in invariant sets of dynamical systems}
\author{T. Medvedev, O. Pochinka}
\date{}
\maketitle
\sloppy

\begin{abstract} The modern qualitative theory of dynamical systems is thoroughly intertwined with the fairly young science of topology. Strange and even bizarre constructions of topology are found sooner or later in dynamics of discrete or continuous dynamical systems. In the present paper we show that the wild Fox-Artin arc naturally emerges in invariant sets of dynamical systems.
\end{abstract}

\section*{Introduction}

Various topological constructions naturally emerge in the modern theory of dynamical systems. For instance, the Cantor set discovered as an example of a set with cardinality of the continuum and zero Lebesgue measure clarified the structure of expanding attractors and contracting repellers. Fractals being the objects of self-similarity with fractional dimension are naturally found in complex dynamics. For example, the basin boundary of an attracting point is  the Julia set.  The lakes of Wada showing the phenomenon of a curve dividing more than two domains were used in construction of the Plykin attractor on the 2-sphere. The irrational winding of the 2-torus being an injectively immersed subset but not a topological submanifold was realized as an invariant manifold of a fixed point of the Anosov diffeomorphism of the 2-torus. The mildly wild embedding of a frame of Debrunner--Fox representing a wild collection of tame arcs in $\mathbb R^3$ was realized as a frame of 1-dimensional saddle points separatrices of Morse-Smale diffeomorphism on the 3-sphere.

\begin{figure}[h]\begin{center}
\epsfig{file=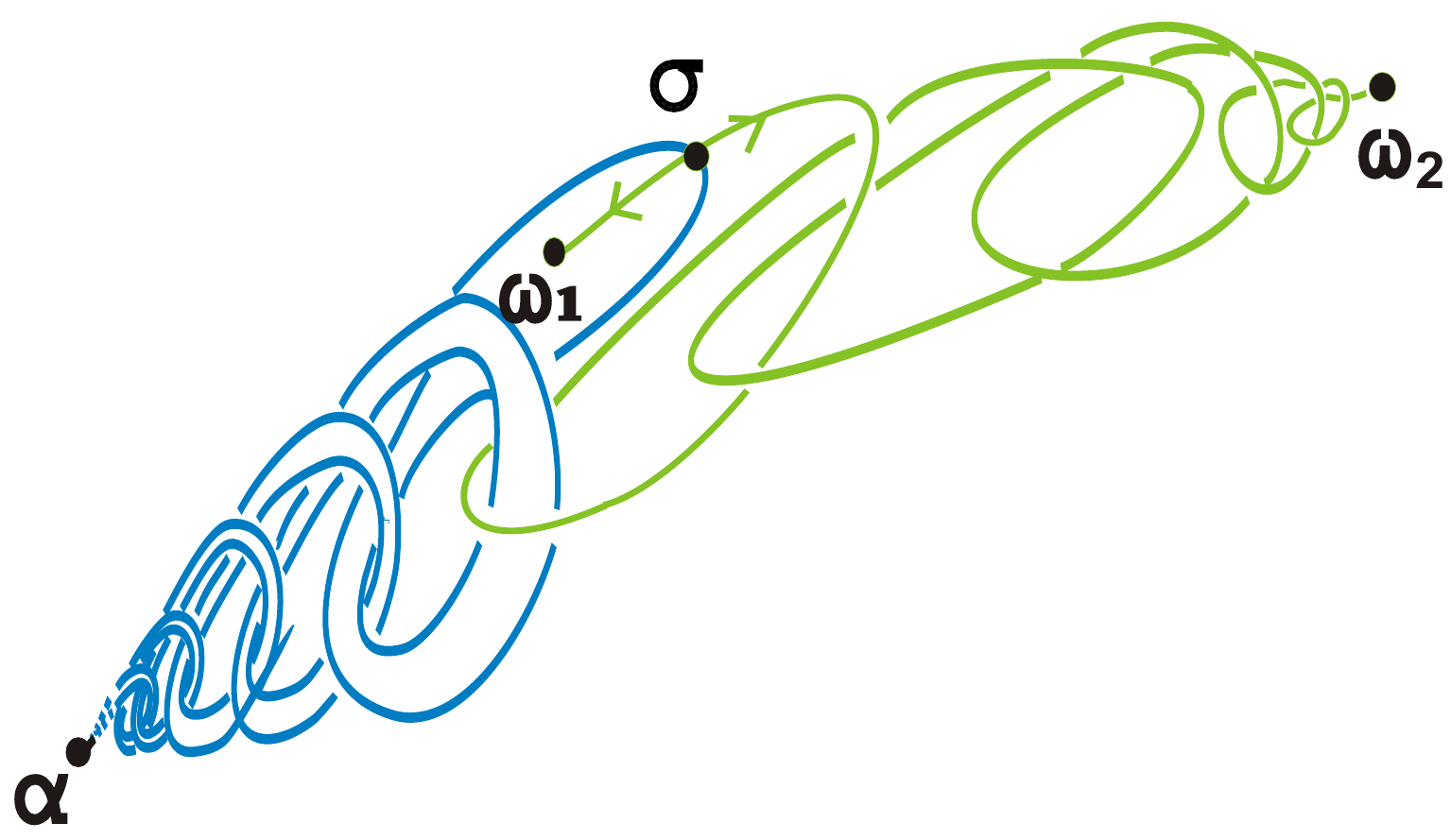, width=10 true cm, height=5. true cm}
\caption{Pixton diffeomorphisn which does not possesses an energy Morse function}\label{ex4}\end{center}
\end{figure}

In the present paper we show the role the wild arc of Fox--Artin \cite{ArFo} plays in dynamics which is yet another example of interconnection between topology and dynamics. For the first time, the Fox--Artin arc appeared in dynamics in 1977 thanks to D. Pixton \cite{Pi}. Pixton constructed a Morse-Smale diffeomorphism  on $3$-sphere with unique saddle point whose one of the unstable separatrices and the stable separatrix form the Fox--Artin arc (see Figure \ref{ex4}). The latter property turns out to be an obstruction to the existence of an energy Morse function for Morse-Smale diffeomorphisms (notice that for any flow an energy function always exists). In 2005 K. Kuperberg \cite{Ku2005} constructed on any orientable 3-manifold without boundary a continuous flow with a discrete set of fixed points and such that the closure of every non-trivial semi-trajectory is the Fox--Artin arc.

In the present paper we show how the Fox--Artin arc naturally emerges as an element of a heteroclinic intersection for regular 4-diffeomorphisms. 

{\it Acknowledgements.} This work was supported by the Russian Foundation for Basic Research (project 16-51-10005-Ko\_a) and the Basic Research Program at the HSE (project 90) in 2017.

\section{The construction of the wild arc of Fox--Artin}

\begin{figure}[ht]
\centerline{\includegraphics[width=11cm,height=7cm]{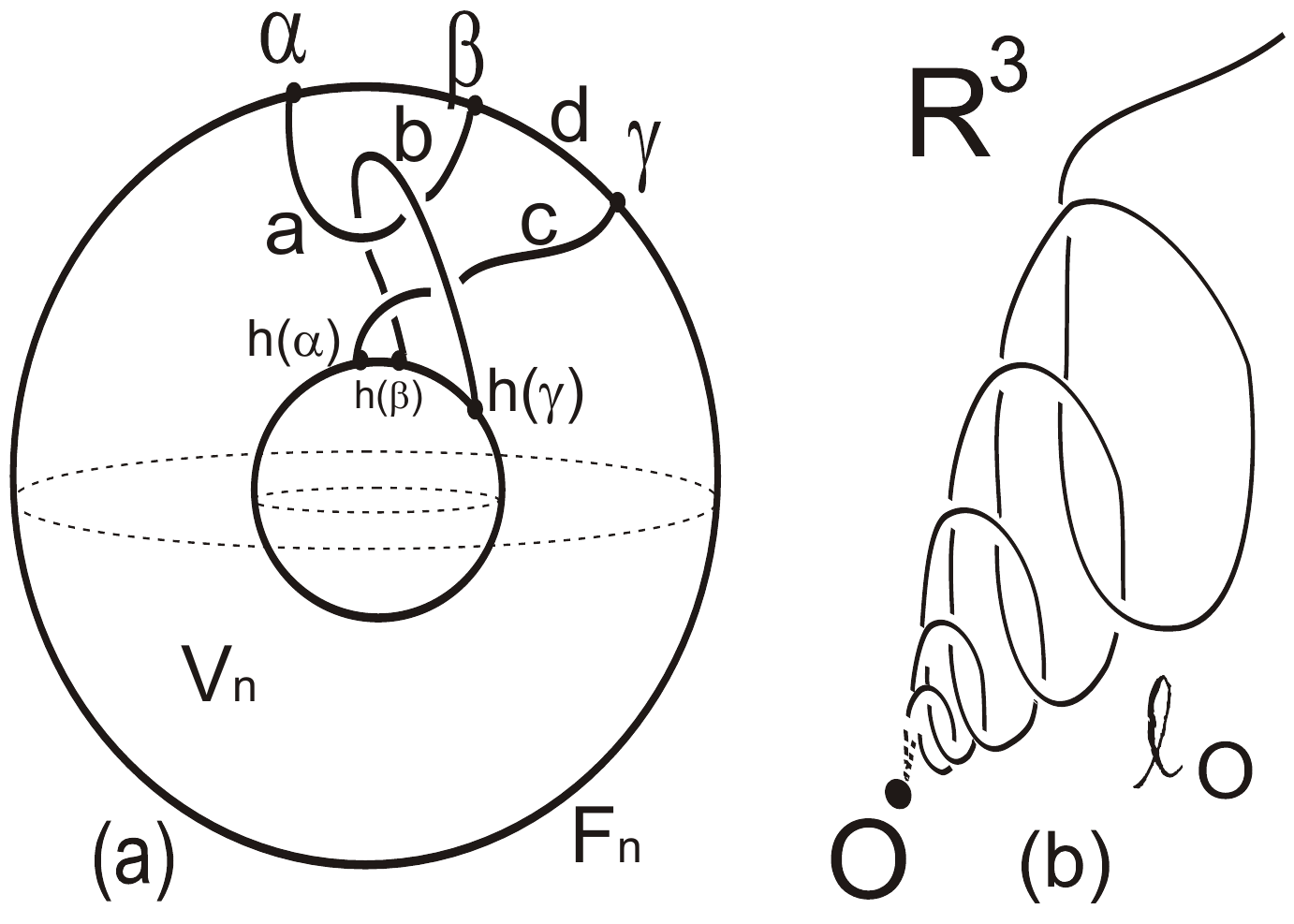}}\caption{Construction of wild arcs in $\mathbb R^3$} \label{wildR3}
\end{figure}

Consider in $\mathbb R^n,n\geq 2$ the $n$-annulus $$V_n=\left\{(x_1,\dots,x_n)\in\mathbb R^n:\frac14\leq x_1^2+\dots+x_n^2\leq 1\right\}$$ and the homothety $h:\mathbb R^n\to\mathbb R^n$ defined by $$h(x_1,\dots,x_n)=\left(\frac{x_1}{2},\dots,\frac{x_n}{2}\right).$$ Let $$F_n=\left\{(x_1,\dots,x_n)\in V_n: x_1^2+\dots+x_n^2= 1\right\},$$ then $\partial V_n=F_n\cup h(F_n)$.
Consider the manifold $\mathbb S^{n-1}\times\mathbb S^1$ as the annulus $V_n$ whose boundaries are identified by $h$. Then the natural projection $$p:\mathbb R^n\setminus O\to\mathbb S^{n-1}\times\mathbb S^1,$$ is a cover. Everywhere below we assume $\mathbb R^k=\{(x_1,\dots,x_n)\in\mathbb R^n:x_{k+1}=\dots=x_n=0\}$ for $k<n$. Further constructions are shown in Figure \ref{wildR3} (a). 

Pick on the circle $F_2$ three different points $\alpha,\beta,\gamma$. Let $d$ be the arc of $F_2$ bounded by $\beta$, $\gamma$ and such that  $\alpha\notin d$. Let $a;b;c\subset V_3$ be pairwise disjoint smooth closed arcs with respective boundary points $\alpha,\beta;h(\beta),h(\gamma);\gamma,h(\alpha)$ such that:
\begin{enumerate}
\item $a,c\subset V_2$, $int~(a\cup c)\cap (F_2\cup h(F_2))=\emptyset$;

\item the closed curve  $b\cup h(d)$ is the boundary of a 2-disk $\Delta\subset V_3$ such that $int~\Delta\cap (F_3\cup h(F_3))=\emptyset$ and $int~\Delta$ intersects both curves $a$ and $c$ only once (one can pick $\Delta$ on the helicoid produced by rotating the curve $h(d)$ around its normal which intersects both $a$ and $c$ at a single point);

\item every connected component of the set $a\cup b\cup c\cup h(a\cup b\cup c)$ is a smooth arc.
\end{enumerate}

Let $\ell=\bigcup\limits_{k\in\mathbb Z}h^k(a\cup b\cup c)$ (see Figure \ref{wildR3} (b)). The arc $\ell_O=\ell\cup O$ is smooth in $\mathbb R^3$ except the point $O$. E.~Artin and R.~Fox proved that $\ell_O$ is wildly embedded, i.e. it is not a submanifold in $\mathbb R^3$ \cite{ArFo}. 

Let $$\hat\ell=p(\ell).$$ One can easily see that $\hat\ell$ is a smooth closed curve (a knot) in $\mathbb S^2\times\mathbb S^1$. This knot was first constructed by B.~Mazur in 1961 \cite{Mazur} as an example of a knot which in $\mathbb S^2\times\mathbb S^1$ is homotopic but not homeomorphic to the trivial knot $\hat l=\{x\}\times\mathbb S^1,x\in\mathbb S^2$. 

Now we show that the constructed above Fox--Artin arc can be realized as invariant subsets of diffeomorphisms with a regular dynamics.

\section{The Fox--Artin arc as an invariant manifold of a saddle point}

The Fox--Artin arc appeared in dynamics for the first time in the paper by D.~Pixton  \cite{Pi} where he constructed a Morse-Smale diffeomorphism of the 3-sphere with a unique saddle point whose invariant manifolds  form the Fox--Artin arc. In his paper \cite{Pi} D.~Pixton showed that because of this ``wild'' dynamics there was no Morse energy function. Later Ch.~Bonatti and V.~Grines classified ``Pixton diffeomorphisms'' in \cite{BoGr} (see also \cite{grin}) and they gave a modern construction of such diffeomorphisms. Now we generalize this construction for arbitrary dimension. It will serve as the base for our examples.

Let $\hat\gamma\subset (\mathbb S^{n-1}\times\mathbb S^1)$ be a homotopy nontrivial smooth knot and let $N(\hat\gamma)$ be its smooth tubular neighborhood. The set $\gamma=p^{-1}(\hat\gamma)$ in $\mathbb R^n$  is an $h$-invariant arc and $N(\gamma)=p^{-1}(N(\hat\gamma))$ is its $h$-invariant tubular neighborhood diffeomorphic to $\mathbb{D}^{n-1}\times\mathbb R^1$. Let {$C=\{(x_2,\dots,x_n)\in\mathbb R^n~:~x_2^2+\dots+x_{n}^2\leq 4\}$} and let the flow $g^t:C\to C$ be defined by $$g^t(x_1,\dots,x_n)=(x_1+t,x_2,\dots,x_n).$$ Then there is a diffeomorphism ${\zeta}:{N(\gamma)}\to C$ that conjugates $h\vert_{{N(\gamma)}}$ and $g=g^1|_C$. 

\begin{figure}[ht]
\centerline{\includegraphics[width=10 true cm, height=7 true cm]{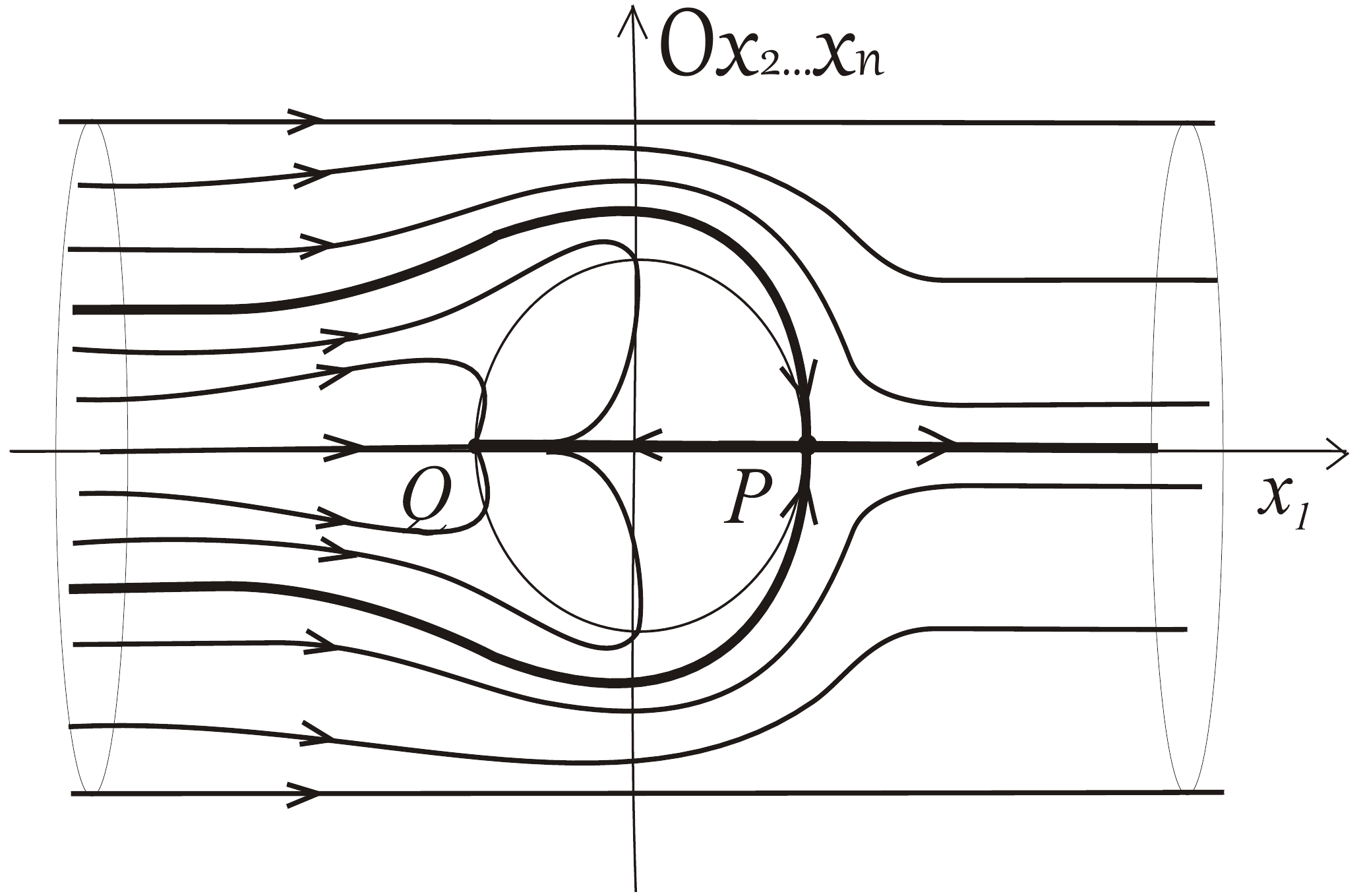}}\caption{Trajectories of the flow $\phi^t$}\label{cherry}
\end{figure}

{Define the flow $\phi^t$ on $C$ by 
$$\begin{cases}
\dot{x}_1=\begin{cases}1-\frac{1}{9}(x_1^2+\dots+x_n^2-4)^2, \quad x_1^2+\dots+x_n^2 \leq 4 \cr
1, \quad otherwise
\end{cases}\cr
\dot{x}_2=\begin{cases}
\frac{x_2}{2}\big(\sin\big(\frac{\pi}{2}\big(x_1^2+\dots+x_n^2-3\big)\big)-1\big), \quad 2<x_1^2+\dots+x_n^2\leq 4\cr
-x_2,\quad \quad x_1^2+\dots+x_n^2\leq 2\cr
0, \quad otherwise
\end{cases}\cr
\dots\dots\dots\dots\dots\dots\cr 
\dot{x}_n=\begin{cases}
\frac{x_n}{2}\big(\sin\big(\frac{\pi}{2}\big(x_1^2+\dots+x_n^2-3\big)\big)-1\big), \quad 2<x_1^2+\dots+x_n^2\leq 4\cr
-x_n,\quad \quad x_1^2+\dots+x_n^2\leq 2\cr
0, \quad otherwise
\end{cases}
\end{cases}$$} 

{By construction the diffeomorphism $\phi=\phi^1$ has exactly two fixed points: the sink $P(1,0,\dots,0)$ and the saddle $Q(-1,0,\dots,0)$ (see Figure \ref{cherry}), both hyperbolic. One unstable separatrix of $Q$ is the open interval $(-1,1)\times\{(0,\dots,0)\}$ belonging to the basin of $P$ and the other is the half-line $(-\infty,-1)\times\{(0,\dots,0)\}$. Notice that $\phi$ coincides with  $g$ outside the ball $\{(x_1,\dots,x_n)\in C:x_1^2+\dots+x_n^2\leq 4\}$}. Define a diffeomorphism  $\bar \varphi:\mathbb R^n\to\mathbb R^n$ in such a way that $\bar{\varphi}$ coincides with $h$ outside ${N(\gamma)}$ and it coincides with  ${\zeta}^{-1}\phi{\zeta}$ on ${N(\gamma)}$. Then $\bar \varphi$ has in ${N(\gamma)}$ exactly two fixed points: the sink ${\zeta}^{-1}(P)$ and the saddle ${\zeta}^{-1}(Q)$, both hyperbolic. The unstable separatrices of ${\zeta}^{-1}(Q)$ lie in $\gamma$.

Denote by $N,S$ the North Pole and the South Pole of the sphere $\mathbb S^n$, respectively. Also denote by 
$\vartheta_S:\mathbb S^n\setminus N\to\mathbb R^n$ and $\vartheta_N:\mathbb S^n\setminus S\to\mathbb R^n$ by the stereographic projections. The diffeomophism $\bar{\varphi}$ induces on $\mathbb{S}^n$ the Morse-Smale diffeomorphism $$\varphi_{\gamma,n}(x)=\begin{cases}\vartheta_S^{-1}(\bar\varphi(\vartheta_S(x))),~x\neq N;\cr N,~x=N\end{cases}.$$ It is immediate from the construction that the non-wandering set of $\varphi_{\gamma,n}$ is four fixed hyperbolic points: two sinks $\omega_{\gamma,n}=\vartheta_S^{-1}({\zeta}^{-1}(P))$, $S_{\gamma,n}=S$, one saddle $\sigma_{\gamma,n}=\vartheta_S^{-1}({\zeta}^{-1}(Q))$ and one source $N_{\gamma,n}=N$.  

Let us choose  natural numbers $k_S$ and $k_N$ such that $\vartheta_S^{-1}(h^{k_S}(V_n))\cap W^s_{\sigma_{\gamma,n}}=\emptyset$ and  $\vartheta_N^{-1}(h^{k_N}(V_n))\cap W^u_{\sigma_{\gamma,n}}=\emptyset$. 
Let $$V_{S_{\gamma,n}}=\vartheta_S^{-1}(h^{k_S}(V_n)),~~~V_{N_{\gamma,n}}=\vartheta_N^{-1}(h^{k_N}(V_n)).$$ Denote by $D_{S_{\gamma,n}}$ the neighborhood of the sink $S_{\gamma,n}$, that is the connected component of $S^n\setminus V_{S_{\gamma,n}}$ and denote by $D_{N_{\gamma,n}}$ the neighborhood of the source $N_{\gamma,n}$, that is the connected component of $S^n\setminus V_{N_{\gamma,n}}$. Then the diffeomorphism $\varphi_{\gamma,n}$ is smoothly conjugated to $h$ in $D_{S_{\gamma,n}}$ by the stereographic projection $\vartheta_S$ and it is smoothly conjugated to $h^{-1}$ in the neighborhood $D_{N_{\gamma,n}}$ of the source $N_{\gamma,n}$ by the stereographic projection $\vartheta_N$.

The standard Pixton diffeomorphism is defined by the map $\varphi_{\ell,3}:\mathbb S^3\to\mathbb S^3$ where $\ell$ is the Fox-Artin arc. Its phase portrait is shown in Figure \ref{ex4}.

\section{Fox-Artin arc as a heteroclinic curve}

In this section we construct an$\Omega$-stable diffeomorphism $f$ on the 4-sphere $S^4$ whose non-wandering set  consists of hyperbolic points: three sinks $\omega_1,\omega_2,\omega_3$, one source $\alpha_1$ and two saddles $\sigma_1,\sigma_2$ such that $dim~W^s_{\sigma_1}=dim~W^s_{\sigma_2}=3$. The wandering set of $f$ contains a heteroclinic intersection $\ell=W^s_{\sigma_1}\cap W^u_{\sigma_2}$ which is the wild Fox-Artin arc in the stable manifold $W^s_{\sigma_1}$ diffeomorphic to $\mathbb R^3$.

Consider the diffeomorphism $\varphi_{l,4}:S^4\to S^4$, where $l=p^{-1}(\hat l)$ and $\hat l$ is a trivial knot in $\mathbb S^2\times\mathbb S^1$ and consider the diffeomorphism $\varphi_{\ell,4}:S^4\to S^4$, where $\ell$ is the Fox-Artin arc.  
It follows from the construction that there exists a diffeomorphism $q:V_{N_{l,4}}\to V_{S_{\ell,4}}$ such that $\vartheta_S(q(W^s_{\sigma_{l,4}}\cap V_{N_{l,4}}))=h^{k_S}(V_3)$ and $q\varphi_{l,4}=\varphi_{\ell,4}q$. Consider  $S^4$ as the union $(S^4\setminus D_{N_{l,4}})\cup_q(S^4\setminus D_{S_{\ell,4}})$ and denote by $\varrho:(S^4\setminus D_{N_{l,4}})\cup(S^4\setminus D_{S_{\ell,4}})\to S^4$ the natural projection. Then the desired difeomorphism $f:S^4\to S^4$ is defined by $$f(x)=\begin{cases}\varrho(\varphi_{l,4}(\varrho^{-1}(x))),~x\in \varrho(S^4\setminus D_{N_{l,4}});\cr \varrho(\varphi_{\ell,4}(\varrho^{-1}(x))),~x\in \varrho(S^4\setminus D_{S_{\ell,4}}).\end{cases}$$

\section{The boundary of a tubular neighborhood of the Fox-Artin arc as a heteroclinic cylinder}

In this section we construct a Morse-Smale diffeomorphism $f$ on the 4-sphere $S^4$ with saddles of codimension one. The non-wandering set $\Omega_f$ consists of two sinks $\omega_1,\omega_2$, two sources $\alpha_1,\alpha_2$ and two saddles  $\sigma_1,\sigma_2$ such that $dim~W^s_{\sigma_1}=dim~W^u_{\sigma_2}=3$. The wandering set of $f$ contains the heteroclinic intersection $H=W^s_{\sigma_1}\cap W^u_{\sigma_2}$ which is a single connected component diffeomorphic to $\mathbb S^1\times\mathbb R$. The cylinder $H$ is embedded into the stable manifold $W^s_{\sigma_1}$, diffeomorphic to $\mathbb R^3$ in such a way that $H$ is the boundary of a tubular neighborhood of a wild Fox-Artin arc in $W^s_{\sigma_1}$.   

Consider the diffeomorphism $\varphi_{l,4}:S^4\to S^4$, where $l=p^{-1}(\hat l)$ and $\hat l$ is a trivial knot in $\mathbb S^2\times\mathbb S^1$ and consider the diffeomorphism $\varphi_{\ell,4}:S^4\to S^4$, where $\ell$ is the Fox-Artin arc. It follows from the construction that there exists a diffeomorphism $r:V_{N_{l,4}}\to V_{N_{\ell,4}}$ such that $\vartheta_N(r(W^s_{\sigma_{l,4}}\cap V_{N_{l,4}}))=h^{k_N}(V_3)$ and $r\varphi^{-1}_{l,4}=\varphi_{\ell,4}r$. Consider  $S^4$ as the union $(S^4\setminus D_{N_{l,4}})\cup_q(S^4\setminus D_{N_{\ell,4}})$ and denote by $\rho:(S^4\setminus D_{N_{l,4}})\cup(S^4\setminus D_{N_{\ell,4}})\to S^4$ the natural projection. Then the desired diffeomorphism $f:S^4\to S^4$ is defined by $$f(x)=\begin{cases}\rho(\varphi^{-1}_{l,4}(\rho^{-1}(x))),~x\in \rho(S^4\setminus D_{N_{l,4}});\cr \rho(\varphi_{\ell,4}(\rho^{-1}(x))),~x\in \rho(S^4\setminus D_{N_{\ell,4}}).\end{cases}$$

\end{document}